\documentclass[reqno,centertags, 12pt]{amsart}
\usepackage{amsmath,amsthm,amscd,amssymb}
\usepackage{latexsym}
\usepackage{graphicx}

\newcommand{\bbC}{{\mathbb{C}}}
\newcommand{\bbD}{{\mathbb{D}}}

\newcommand{\bbR}{{\mathbb{R}}}

\newcommand{\bbZ}{{\mathbb{Z}}}

\newcommand{\calS}{{\mathcal S}}

\newcommand{\lb}{\label}
\newcommand{\f}{\frac}

\newcommand{\ol}{\overline}
\newcommand{\ti}{\tilde  }

\newcommand{\tr}{\text{\rm{Tr}}}

\newcommand{\supp}{\text{\rm{supp}}}

\newcommand{\bi}{\bibitem}

\newcommand{\beq}{\begin{equation}}
\newcommand{\eeq}{\end{equation}}
\newcommand{\ba}{\begin{align}}
\newcommand{\ea}{\end{align}}
\newcommand{\veps}{\varepsilon}

\DeclareMathOperator{\Real}{Re}

\DeclareMathOperator{\Ima}{Im}

\allowdisplaybreaks
\numberwithin{equation}{section}

\newtheorem{theorem}{Theorem}[section]

\newtheorem{proposition}[theorem]{Proposition}
\newtheorem{lemma}[theorem]{Lemma}

\theoremstyle{definition}

\newtheorem{example}[theorem]{Example}

\theoremstyle{remark}

\newcommand{\abs}[1]{\lvert#1\rvert}

\newcounter{smalllist}
\newenvironment{SL}{\begin{list}{{\rm\roman{smalllist})}}{%
\setlength{\topsep}{0mm}\setlength{\parsep}{0mm}\setlength{\itemsep}{0mm}%
\setlength{\labelwidth}{2em}\setlength{\leftmargin}{2em}\usecounter{smalllist}%
}}{\end{list}}

\begin{document}

\title[Fine Structure of the Zeros of OP, III]
{Fine Structure of the Zeros of Orthogonal Polynomials, \\
III. Periodic Recursion Coefficients}
\author[B.~Simon]{Barry Simon*}

\thanks{$^*$ Mathematics 253-37, California Institute of Technology, Pasadena, CA 91125, USA. 
E-mail: bsimon@caltech.edu. Supported in part by NSF grant DMS-0140592} 
\thanks{To be submitted to Comm. Pure Appl. Math.}

\date{December 16, 2004}

\begin{abstract} We discuss asymptotics of the zeros of orthogonal polynomials on the 
real line and on the unit circle when the recursion coefficients are periodic. The zeros 
on or near the absolutely continuous spectrum have a clock structure with spacings 
inverse to the density of zeros. Zeros away from the a.c.\ spectrum have limit points 
mod $p$ and only finitely many of them.  
\end{abstract}

\maketitle

\section{Introduction} \lb{s1} 

This paper is the third in a series \cite{SaffProc,Saff2} that discusses detailed 
asymptotics of the zeros of orthogonal polynomials with special emphasis on 
distances between nearby zeros. We discuss both orthogonal polynomials on the 
real line (OPRL) where the basic recursion for the orthonormal polynomials, 
$p_n(x)$, is 
\begin{equation} \lb{1.1} 
xp_n(x) = a_{n+1} p_{n+1}(x) + b_{n+1} p_n (x) + a_n p_{n-1}(x) 
\end{equation} 
($a_n >0$ for $n=1,2, \dots$, $b_n$ real, and $p_{-1}(x)\equiv 0$), and 
orthogonal polynomials on the unit circle (OPUC) where the basic recursion is 
\begin{equation} \lb{1.2} 
\varphi_{n+1}(z) = \rho_n^{-1} (z\varphi_n (z) -\bar\alpha_n \varphi_n^*(z)) 
\end{equation} 
Here $\alpha_n$ are complex coefficients lying in the unit disk $\bbD$ and 
\begin{equation} \lb{1.3} 
\varphi_n^*(z) = z^n \, \ol{\varphi_n (1/\bar z)} 
\end{equation} 
and 
\begin{equation} \lb{1.4} 
\rho_n =(1-\abs{\alpha_n}^2)^{1/2} 
\end{equation} 

In this paper, we focus on the case where the Jacobi coefficients $\{a_n\}_{n=1}^\infty, 
\{b_n\}_{n=1}^\infty$ or the Verblunsky coefficients $\{\alpha_n\}_{n=0}^\infty$ are 
periodic, that is, for some $p$, 
\begin{equation} \lb{1.5} 
a_{n+p} =a_n \qquad b_{n+p} =b_n 
\end{equation} 
or 
\begin{equation} \lb{1.6} 
\alpha_{n+p}=\alpha_n 
\end{equation}
It should be possible to say something about perturbations of a periodic sequence, 
say $\alpha_n^{(0)}$, which obeys \eqref{1.6} and $\alpha_n = \alpha_n^{(0)} + 
\delta\alpha_n$ with $\abs{\delta\alpha_n}\to 0$ sufficiently fast. We leave 
the details to be worked out elsewhere. 

To describe our results, we begin by summarizing some of the basics of the structure of 
the measures and recursion relations when \eqref{1.5} or \eqref{1.6} holds. We will 
say more about this underlying structure in the sections below. In this introduction, 
we will assume that all gaps are open, although we don't need and won't use that 
assumption in the detailed discussion. 

When \eqref{1.5} holds, the continuous part of the underlying measure, $d\rho$, on 
$\bbR$ is supported on $p$ closed intervals $[\alpha_j, \beta_j]$, $j=1, \dots, p$,  
called bands, with gaps $(\beta_j, \alpha_{j+1})$ in between. Each gap has zero or one  
mass point. The $m$-function of the measure $d\rho$, 
\begin{equation} \lb{1.7} 
m(z) = \int \f{d\rho(x)}{x-z} 
\end{equation} 
has a meromorphic continuation to the genus $p-1$ hyperelliptic Riemann surface, $\calS$, 
associated to $[\prod_{j=1}^p (x-\alpha_j)(x-\beta_j)]^{1/2}$. This surface has a 
natural projection $\pi :\calS\to\bbC$, a twofold cover except at the branch points 
$\{\alpha_j\}_{j=1}^p \cup \{\beta_j\}_{j=1}^p$. $\pi^{-1} [\beta_j, \alpha_{j+1}]$ 
is a circle and $m(z)$ has exactly one pole $\gamma_1, \dots, \gamma_{p-1}$ on each 
circle. 

It has been known for many years (see Faber \cite{Fab}) that 
the density of zeros $dk$ is supported on $\cup_{j=1}^p [\alpha_j, \beta_j]\equiv B$ 
and is the equilibrium measure for $B$ in potential theory. We define $k(E) = 
\int_{\alpha_1}^E dk$. Then $k(\beta_j)= j/p$. Our main results about OPRL are: 
\begin{SL} 
\item[(1)] We can describe the zeros of $p_{np-1}(x)$ exactly (not just asymptotically) 
in terms of $\pi (\gamma_j)$ and $k(E)$. 

\item[(2)] Asymptotically, as $n\to\infty$, the number of zeros of $p_n$ in each band 
$[\alpha_j, \beta_j]$, $N^{(n,j)}$, obeys $\sup_n \abs{\f{n}{p} - N^{(n,j)}}<\infty$,  
and the zeros $\{x_\ell^{(n,j)}\}_{\ell=1}^{N(n,j)}$ obey 
\begin{equation} \lb{1.8} 
\sup_{\substack{ j \\ \ell=1,2, \dots, N^{(n,j)}-1}} \, n\biggl| k(x_{\ell+1}^{(n,j)}) - 
k(x_\ell^{(n,j)}) - \f{1}{n}\biggr| \to 0
\end{equation} 
as $n\to\infty$.  

\item[(3)] $z\in\bbC$ is a limit of zeros of $p_n$ if and only if $z$ lies in 
$\supp (d\rho)$.

\item[(4)] Outside the bands, there are at most $2p+2b-3$ points which are limits of 
zeros of $p_{mp+b-1}$ for each $b=1, \dots, p$ and, except for these limits, zeros 
have no accumulation points in $\bbC\backslash\text{bands}$. 
\end{SL} 

\smallskip
For OPUC, the continuous part of the measure, $d\mu$, is supported  on $p$ disjoint 
intervals $\{e^{i\theta}\mid x_j \leq \theta \leq y_j\}$, $j=1, \dots, p$, in $\partial\bbD$ 
with $p$ gaps in between $\{e^{i\theta}\mid y_j\leq \theta \leq x_{j+1}\}$ with 
$x_{p+1}\equiv 2\pi + x_1$. Each gap has zero or one mass point. The Carath\'eodory 
function of the measure $d\mu$, 
\begin{equation} \lb{1.10} 
F(z) =\int \f{e^{i\theta}+z}{e^{i\theta}-z}\, d\mu(\theta) 
\end{equation} 
has a meromorphic continuation from $\bbD$ to the genus $p-1$ hyperelliptic 
Riemann surface, $\calS$, associated to $[\prod_{j=1}^p (z-e^{ix_j})(z-e^{iy_j})]^{1/2}$. 
The surface has a natural projection $\pi:\calS\to\bbC$, and the closure of each gap 
has a circle as the inverse image. $F$ has a single pole in each such circle, so $p$ 
in all at $\gamma_1, \dots, \gamma_p$. 

Again, the density of zeros is the equilibrium measure for the bands and each band 
has mass $1/p$ in this measure. See \cite{OPUC2}, especially Chapter~11, for a discussion 
of periodic OPUC. Our main results for OPUC are: 
\begin{SL} 
\item[($1'$)] We can describe the zeros of $\varphi_{np}^* -\varphi_{np}$ exactly 
(note, not zeros of $\varphi_{np}$). 

\item[($2'$)] Asymptotically, as $n\to\infty$, the number of zeros of $\varphi_n$ near 
each band, $N^{(n,j)}$, obeys $\sup_n \abs{\f{n}{p} - N^{(n,j)}}<\infty$, and the points 
on the bands closest to the zeros obey an estimate like \eqref{1.8}.  

\item[($3'$)] $z\in\bbC$ is a limit of zeros of $\varphi_n$ if and only if $z$ lies 
in $\supp (d\mu)$. 

\item[($4'$)] There are at most $2p+2b-1$ points which are limits of zeros of 
$\varphi_{mp+b}$ for each $b=1, \dots, p$ and, except for these limits, zeros have 
no accumulation points in $\bbC\backslash\text{bands}$. 
\end{SL} 

\smallskip
In Section~\ref{s2}, we discuss OPRL when \eqref{1.5} holds, and in Section~\ref{s3}, 
OPUC when \eqref{1.6} holds. Each section begins with a summary of transfer matrix 
techniques for periodic recursion coefficients (Floquet theory). 

While I am unaware of any previous work on the precise subject of Sections~\ref{s2} 
and \ref{s3}, the results are closely related to prior work of Peherstorfer 
\cite{Peh03,Peh04}, who discusses zeros in terms of measures supported on a union 
of bands with a  particular structure that overlaps our class of measures. For a 
discussion of zeros for OPUC with two bands, see \cite{LukPeh}. 

These papers also consider situations where the recursion coefficients are only almost 
periodic. For any finite collection of closed intervals on $\bbR$ or closed arcs on 
$\partial\bbD$, there is a natural isospectral torus of OPRL or OPUC where the 
corresponding $m$- or $F$-function has minimal degree on the Riemann surface 
(see, e.g., \cite[Section~11.8]{OPUC2}). It would be interesting to extend the 
results of the current paper to that case. 

\smallskip 

It is a pleasure to thank Chuck Newman and Percy Deift for the hospitality of the 
Courant Institute where some of this work was done. 

\section{OPRL With Periodic Jacobi Coefficients} \lb{s2} 

In this section, we analyze the zeros of OPRL with Jacobi coefficients obeying \eqref{1.5}. 
We begin with a summary of the theory of transfer matrices, discriminants, and Abelian 
functions associated to this situation. A reference for much of this theory is 
von Moerbeke \cite{vMoer}; a discussion of the discriminant can be found in 
Hochstadt \cite{Hoch7}, von Moerbeke \cite{vMoer}, Toda \cite{Toda}, and Last 
\cite{Last92}. The theory is close to the OPUC theory developed in Chapter~11 
of \cite{OPUC2}.  

Define the $2\times 2$ matrix,  
\begin{equation} \lb{2.1} 
A_k(z) = \f{1}{a_{k+1}}\, \begin{pmatrix} 
z-b_{k+1} & -a_k \\
a_{k+1} & 0 \end{pmatrix} 
\end{equation} 
where 
\begin{equation} \lb{2.2} 
a_0 \equiv a_p 
\end{equation} 
Thus 
\begin{equation} \lb{2.3} 
\det (A_k) = \f{a_k}{a_{k+1}} 
\end{equation} 
and the abstract form of \eqref{1.1} 
\begin{equation} \lb{2.4} 
zu_n = a_{n+1} u_{n+1} + b_{n+1} u_n + a_n u_{n-1}  
\end{equation} 
is equivalent to 
\begin{equation} \lb{2.5} 
\binom{u_{n+1}}{u_n} = A_n \binom{u_n}{u_{n-1}}  
\end{equation} 
So, in particular, 
\begin{equation} \lb{2.6} 
\binom{p_{n+1}(z)}{p_n(z)} = A_n A_{n-1} \dots A_0 \binom{1}{0} 
\end{equation}

This motivates the definition of the transfer matrix,  
\begin{equation} \lb{2.7} 
T_n(z) = A_{n-1}(z) \dots A_0 (z)  
\end{equation}
for $n=1,2,\dots$. We have, by \eqref{1.5}, that 
\begin{equation} \lb{2.8} 
T_{mp+b} = T_b (T_p)^m 
\end{equation} 
suggesting that $T_p$ plays a basic role. By \eqref{2.3} and \eqref{2.2}, 
\begin{equation} \lb{2.9} 
\det (T_p) =1 
\end{equation} 

A fundamental quantity is the discriminant 
\begin{equation} \lb{2.10} 
\Delta(z) = \tr (T_p(z)) 
\end{equation} 
By \eqref{2.6}, we have 
\begin{equation} \lb{2.10a} 
T_n(z) = \begin{pmatrix} p_n(z) & q_{n-1}(z) \\
p_{n-1}(z) & q_{n-2}(z) \end{pmatrix} 
\end{equation} 
where $q_n(z)$ is a polynomial of degree $n$ that is essentially the polynomial 
of the second kind (the normalization is not the standard one but involves an 
extra $a_p$). 

By \eqref{2.9} and \eqref{2.10}, $T_p(z)$ has eigenvalues 
\begin{equation} \lb{2.11} 
\Gamma_\pm (z) = \f{\Delta(z)}{2} \pm \sqrt{\bigl(\tfrac{\Delta (z)}{2}\bigr)^2 -1}
\end{equation} 
In a moment, we will define branch cuts in such a way that on all of $\bbC\backslash 
\text{cuts}$, 
\begin{equation} \lb{2.12x} 
\abs{\Gamma_+(z)} > \abs{\Gamma_- (z)} 
\end{equation} 
so \eqref{2.8} implies the Lyapunov exponent is given by 
\begin{equation} \lb{2.12}
\lim_{n\to\infty} \, \f{1}{n} \, \log \|T_n(z)\| = \f{1}{p} \, \log \abs{\Gamma_+(z)} 
\equiv \gamma(z) 
\end{equation} 

\eqref{2.11} means $\abs{\Gamma_+}=\abs{\Gamma_-}$ if and only if $\Delta(z)\in 
[-2,2]$, and one shows that this only happens if $z$ is real. Moreover, if $\Delta(z)\in 
(-2,2)$, then $\Delta'(x)\neq 0$. Thus, for $x$ very negative, $(-1)^p \Delta(x) >0$ and 
solutions of $(-1)^p \Delta(x)=\pm 2$ alternate as $+2,-2,-2,+2,+2,-2,-2,\dots$, which 
we label as 
\begin{equation} \lb{2.13} 
\alpha_1 < \beta_1 \leq \alpha_2 < \beta_2 \leq\alpha_3 < \cdots <\beta_p
\end{equation} 
Since $\Delta(x)$ is a polynomial of degree $p$, there are $p$ solutions of 
$\Delta(x)=2$ and of $\Delta(x)=-2$, so $2p$ points $\{\alpha_j\}_{j=1}^p \cup 
\{\beta_j\}_{j=1}^p$. 

The bands are $[\alpha_1,\beta_1],[\alpha_2,\beta_2],\dots, [\alpha_p,\beta_p]$ and the 
gaps are $(\beta_1,\alpha_2)$, $(\beta_2,\alpha_3), \dots, (\beta_{p-1}, \alpha_{p+1})$. 
If some $\beta_j=\alpha_{j+1}$, we say the $j$-th gap is closed. Otherwise 
we say the gap is open. 

If we remove the bands from $\bbC$, $\Gamma_\pm(z)$ are single-valued analytic functions 
and \eqref{2.12x} holds. Moreover, $\Gamma_+$ has an analytic continuation to the Riemann 
surface, $\calS$, of genus $\ell \leq p-1$ where $\ell$ is the number of open gaps. 
$\calS$ is defined by the function $[(z-\alpha_1)(z-\beta_p) \prod_{\text{open gaps}} 
(z-\beta_j) (z-\alpha_{j+1})]^{1/2}$. $\Gamma_-$ is precisely the analytic continuation 
of $\Gamma_+$ to the second sheet. 

The Dirichlet data are partially those $x$'s where 
\begin{equation} \lb{2.14} 
T_p(x) \binom{1}{0} = c_x \binom{1}{0}  
\end{equation} 
that is, points where the $21$ matrix element of $T_p$ vanishes. It can be seen that 
the Dirichlet data $x$'s occur, one to each gap, that is, $x_1, \dots, x_{p-1}$ with 
$\beta_j \leq x_j \leq \alpha_{j+1}$. If $x$ is at an edge of a gap, then $c_j \equiv 
c_{x_j}$ is $\pm 1$. Otherwise $\abs{c_j} \neq 1$. If $\abs{c_j} >1$, we add the 
sign $\sigma_j =-1$ to $x_j$, and if $\abs{c_j} <1$, we add the sign $\sigma_j= +1$ 
to $x_j$. Thus the values of Dirichlet data for each open gap are two copies of $[\alpha_j, 
\beta_j]$ glued at the ends, that is, a circle. The set of Dirichlet data is thus 
an $\ell$-dimensional torus. It is a fundamental result \cite{vMoer} that the map from 
$a$'s and $b$'s to Dirichlet data sets up a one-one correspondence to all $a$'s and 
$b$'s with a given $\Delta$, that is, the set of $a$'s and $b$'s with a given $\Delta$ 
is an $\ell$-dimensional torus. 

The $m$-function \eqref{1.7} associated to $d\rho$ has a meromorphic continuation to 
the Riemann surface, $\calS$, with poles precisely at the points $x_j$ on the 
principal sheet if $\sigma_j=+1$ and on the bottom sheet if $\sigma_j =-1$. $\rho$ has 
point mass precisely at those $x_j\in (\beta_j, \alpha_{j+1})$ with $\sigma_j =+1$. It 
has absolutely continuous support exactly the union of the bands, and has no singular 
part other than the possible point masses in the gaps. 

Finally, in the review, we note that the potential theoretic equilibrium measure $dk$ 
for the set of bands has several critical properties:  
\begin{SL} 
\item[(1)] If $k(x) = \int_{\alpha_1}^x dk$, then 
\begin{equation} \lb{2.15} 
k(\beta_j) = k(\alpha_{j+1}) = \f{j}{p} 
\end{equation} 
\item[(2)] The Thouless formula holds: 
\begin{equation} \lb{2.16} 
\gamma(z) =\int \log\abs{z-x}\, dk(x) + \log C_B 
\end{equation} 
where $\gamma$ is given by \eqref{2.12} and $C_B$ is the (logarithmic) 
capacity of $B$. 
\item[(3)] The (logarithmic) capacity of the bands is given by 
\begin{equation} \lb{2.17} 
C_B = \biggl(\, \prod_{j=1}^p a_j\biggr)^{-1} 
\end{equation} 
\item[(4)] 
\begin{equation} \lb{2.18x} 
\Gamma_+(z) =C_B \exp \biggl( p \int \log (z-x)\, dk(x) \biggr)
\end{equation} 
\end{SL} 

That completes the review of periodic OPRL. We now turn to the study of the 
zeros. We begin by describing exactly (not just asymptotically!) the zeros of 
$P_{mp-1}$: 

\begin{theorem}\lb{T2.1} The zeros of $P_{mp-1}(x)$ are exactly 
\begin{SL} 
\item[{\rm{(i)}}] The $p-1$ Dirichlet data points $\{x_j\}_{j=1}^{p-1}$. 
\item[{\rm{(ii)}}] The $(m-1)p$ points $\{x_{k,q}^{(m)}
\}_{\substack{k=1, \dots, p \\ q=1, \dots, m-1}}$ where 
\begin{equation} \lb{2.18} 
k(x_{k,q}^{(m)}) = \f{k-1}{p} + \f{q}{mp} 
\end{equation} 
\end{SL} 
\end{theorem} 

{\it Remarks.} 1. The points of \eqref{2.18} can be described as follows. Break 
each band $[\alpha_j,\beta_j]$ into $m$ pieces of equal size in equilibrium measure. 
The $x_{k,q}^{(m)}$ are the interior break points. 

\smallskip 
2. If a gap is closed, we include its position in the ``Dirichlet points" of (i). 

\smallskip 
3. Generically, there are not zeros at the band edges, that is, \eqref{2.18} has $q=1, 
\dots, m-1$ but not $q=0$ or $q=m$. But it can happen that one or more of the 
Dirichlet data points is at an $\alpha_{j+1}$ or a $\beta_j$. 

\smallskip 
4. This immediately implies that once one proves that the density of zeros exists, 
that it is given by $dk$. 

\smallskip 
5. It is remarkable that this result is new, given that it is so elegant and its 
proof so simple! I think this is because the OP community most often focuses on 
measures and doesn't think so much about the recursion parameters and the Schr\"odinger 
operator community doesn't usually think of zeros of $P_n$. 

\begin{example} \lb{E2.2} Let $b_n\equiv 0$, $a_n \equiv \f12$ which has period 
$p=1$. It is well-known in this case that the $P_n$ are essentially Chebyshev 
polynomials of the second kind, that is, 
\begin{equation} \lb{2.19} 
P_n (\cos\theta) = \f{1}{2^n}\, \f{\sin(n+1)\theta}{\sin\theta} 
\end{equation} 
Thus $P_{m-1}$ has zeros at points where 
\begin{equation} \lb{2.20x} 
\theta =\f{j\pi}{m} \qquad j=1, \dots, m-1 
\end{equation} 
(the zeros at $\theta=0$ and $\theta=\pi$ are cancelled by the $\sin(\theta)$).   
$k(x)=\pi -\arccos (x)$ and \eqref{2.20x} is \eqref{2.18}. We see that 
Theorem~\ref{T2.1} generalizes the obvious result on the zeros of the 
Chebyshev polynomials of the second kind. \qed 
\end{example} 

\begin{proof}[First Proof of Theorem~\ref{T2.1}] By \eqref{2.10a}, zeros of $P_{mp-1}$ 
are precisely points where the $12$ matrix element of $T_{mp}$ vanishes, that is, 
points where $\binom{1}{0}$ is an eigenvector of $T_{mp}$. That is, zeros of 
$P_{mp-1}$ are Dirichlet points for this period $mp$ problem. 

When \eqref{1.5} holds, we can view the $a$'s and $b$'s as periodic of period $mp$. 
There are closed gaps where $T_{mp}(z)=\pm\boldsymbol{1}$, that is, interior points 
to the original bands where $(\Gamma_\pm)^m =1$, that is, points where \eqref{2.18} 
holds. Thus, the Dirichlet data for $T_{mp}$ are exactly the points claimed. 
\end{proof}

Theorem~\ref{T2.1} immediately implies point (2) from the introduction. 

\begin{theorem} \lb{T2.3} Let $P_n(x)$ be a family of OPRL associated to a 
set of Jacobi parameters obeying \eqref{1.5}. Let $(\alpha_j,\beta_j)$ be 
a single band and let $N^{(n,j)}$ be the number of zeros of $P_n$ in that band. 
Then 
\begin{equation} \lb{2.20} 
\abs{N^{(mp+b,j)} - (m-1)} \leq \min (b+1, p-b) 
\end{equation} 
for $-1\leq b\leq p-1$. In particular, 
\begin{equation} \lb{2.21} 
\biggl| N^{(n,j)} - \f{n}{p}\biggr| \leq 1 + \f{p}{2} 
\end{equation} 
\end{theorem} 

\begin{proof} By a variational principle for any $n,n'$, 
\begin{equation} \lb{2.22} 
\abs{N^{(n,j)} - N^{(n',j)}} \leq \abs{n-n'} 
\end{equation} 
\eqref{2.20} is immediate from Theorem~\ref{T2.1} if we take $n'=mp-1$ and $n'=
mp+(p-1)$. \eqref{2.21} follows from \eqref{2.20} given that $\min (b+1,p-b) 
\leq p/2$.  
\end{proof} 

{\it Remark.} Because of possibilities of Dirichlet data zeros at $\alpha_j$ and/or 
$\beta_j$, we need $(\alpha_j, \beta_j)$ in defining $N^{(n,j)}$. It is more natural 
to use $[\alpha_j, \beta_j]$. If one does that, \eqref{2.20} becomes $2+\min(b+1,p-b)$ 
and \eqref{2.21}, $3+\f{p}{2}$. 

\smallskip 
To go beyond these results and prove clock behavior for the zeros of $p_{mp+b}$ 
($b\not\equiv -1$ mod $p$), we need to analyze the structure of $p_n$ in terms of 
$\Gamma_+, \Gamma_-$. For $z$ not a branch point (or closed gap), $\Gamma_+ \neq 
\Gamma_-$. $\Gamma_+$ is well-defined on $\bbC\backslash\text{bands}$ since  
$\abs{\Gamma_+} > \abs{\Gamma_-}$. On the bands, $\abs{\Gamma_+}=\abs{\Gamma_-}$  
and, indeed, the boundary values on the two sides of a band are distinct. But 
$\Gamma_+$ is analytic on $\bbC\backslash\text{bands}$, so for such $z$, we 
can define $P_\pm$ by 
\begin{equation} \lb{2.23} 
T_p(z) = \Gamma_+ P_+ + \Gamma_- P_- 
\end{equation} 
where $P_+, P_-$ are $2\times 2$ rank one projections obeying 
\begin{equation} \lb{2.24} 
P_+^2 = P_+ \qquad P_-^2 =P_- \qquad P_+ P_- = P_- P_+ =0 
\end{equation} 
and 
\begin{equation} \lb{2.25} 
P_+ + P_- = \boldsymbol{1} 
\end{equation} 
It follows from \eqref{2.23} and \eqref{2.25} that 
\begin{align} 
P_+ &= \f{T_p(z) - \Gamma_- (z) \boldsymbol{1}}{\Gamma_+ - \Gamma_-} \lb{2.26} \\
P_- &= \f{T_p(z) - \Gamma_+ (z) \boldsymbol{1}}{\Gamma_- - \Gamma_+} \lb{2.27} 
\end{align} 
which, in particular, shows that $P_+$ is a meromorphic function on $\calS$ whose 
second-sheet values are just $P_-$. 

Define 
\begin{align} 
a(z) &= \bigg\langle \binom{0}{1}, P_+(z) \binom{1}{0}\bigg\rangle \lb{2.28} \\
b(z) &= \bigg\langle \binom{1}{0}, P_+(z) \binom{1}{0}\bigg\rangle \lb{2.29} 
\end{align} 
so \eqref{2.25} implies 
\begin{align} 
\biggl\langle \binom{0}{1}, P_-(z) \binom{1}{0}\bigg\rangle &= -a(z) \lb{2.30} \\
\bigg\langle \binom{1}{0}, P_-(z) \binom{1}{0}\bigg\rangle &= 1-b(z) \lb{2.31} 
\end{align} 

Under most circumstances, $a(z)$ has a pole at band edges where $\Gamma_+ -\Gamma_- 
\to 0$. For later purpose, we note that $\langle \binom{0}{1}, (T_p(z) - \Gamma_- 
\boldsymbol{1}) \binom{1}{0}\rangle = \langle \binom{0}{1}, T_p(z) \binom{1}{0} 
\rangle$ has a finite limit at such points. Later we will be looking at 
\begin{align*}
a(z) (\Gamma_+^m - \Gamma_-^m) &= \biggl\langle\binom{0}{1}, T_p(z)\binom{1}{0} 
\bigg\rangle \, \f{\Gamma_+^m - \Gamma_-^m}{\Gamma_+ - \Gamma_-} \\ 
&  \to \bigg\langle \binom{0}{1}, T_p(z) \binom{1}{0} \bigg\rangle m\Gamma_+^{m-1} 
\end{align*}
if $\Gamma_+ - \Gamma_- \to 0$. This is zero if and only if $\langle \binom{0}{1}, 
T_p(z) \binom{1}{0}\rangle =0$, that is, if and only if the edge of the band is 
a Dirichlet data point.  

\eqref{2.23} and \eqref{2.24} imply 
\begin{equation} \lb{2.32} 
T_{mp}(z) = T_p(z)^m = \Gamma_+^m P_+ + \Gamma_-^m P_- 
\end{equation} 
so 
\begin{equation} \lb{2.33} 
T_{mp}(z)\binom{0}{1} = [a(z)(\Gamma_+^m - \Gamma_-^m)]\binom{0}{1} + 
[b(z) \Gamma_+^m + (1-b(z)) \Gamma_-^m)] \binom{1}{0}  
\end{equation} 
Thus, by \eqref{2.21} for $b\geq 0$, 
\begin{align} 
P_{mp+b-1} &= \bigg\langle \binom{0}{1}, T_b T_{mp}\binom{1}{0}\bigg\rangle \lb{2.37a} \\
&= [(\Gamma_+^m - \Gamma_-^m)a(z)] q_{b-2}(z) + [b(z) \Gamma_+^m + (1-b(z)) 
\Gamma_-^m)] p_{b-1}(z) \lb{2.34} 
\end{align} 
where 
\begin{equation} \lb{2.35} 
q_{-2}(z) \equiv 1 \qquad q_{-1}(z) \equiv 0 
\end{equation} 

\begin{proof}[Second Proof of Theorem~\ref{T2.1}] For $b=0$, $p_{b-1}\equiv 0$ and 
$q_{b-2}=1$, so 
\begin{equation} \lb{2.36} 
p_{mp-1}(z) = (\Gamma_+^m - \Gamma_-^m) a(z) 
\end{equation} 
Its zeros are thus points where $a(z)=0$ or where $\Gamma_+^m = \Gamma_-^m$, except 
that at branch points, $a(z)$ can have a pole which can cancel a zero of $\Gamma_+^m 
-\Gamma_-^m$. 

$a(z)=0$ if and only if $\binom{1}{0}$ is an eigenvector of $T_p(z)$, that is, 
exactly at the Dirichlet data points. 

$\Gamma_+^m = \Gamma_-^m$ is equivalent to $\Gamma_+^{2m}=1$ since $\Gamma_- = 
\Gamma_+^{-1}$. This implies $\abs{\Gamma_+} = \abs{\Gamma_-}$, so can only happen 
on the bands. On the bands, by \eqref{2.20}, 
\begin{equation} \lb{2.37} 
\Gamma_+(x) = \exp (\pi i p\, k(x)) 
\end{equation} 
and $\Gamma_+^{2m}=1$ if and only if 
\begin{equation} \lb{2.38} 
mp \, k(x) \in \bbZ
\end{equation} 
that is, if \eqref{2.18} holds for some $q=0, \dots, m$. But at $q=0$ or $q=m$, 
$a(z)$ has a pole that cancels the zero of $\Gamma_+^m - \Gamma_-^m$, so the zeros 
of $p_{mp-1}$ are precisely given by (i) and (ii) of Theorem~\ref{T2.1}. 
\end{proof} 

We can use \eqref{2.34} to analyze zeros of $p_{mp+b-1}$ for large $m$. We begin 
with the region away from the bands: 

\begin{theorem} \lb{T2.4} Let $z\in\bbC\backslash\text{bands}$ and let $b$ be fixed. 
Then 
\begin{equation} \lb{2.39} 
\lim_{m\to\infty} \, \Gamma_+(z)^{-m} p_{mp+b-1}(z) = a(z) q_{b-2}(z) 
+ b(z) p_{b-1}(z) 
\end{equation} 
In particular, if the right side of \eqref{2.39} is called $j_b(z)$, then  
\begin{SL}  
\item[{\rm{(1)}}] If $j_b(z_0)\neq 0$, then $p_{mp+b-1}(z)$ is nonvanishing near $z_0$ 
for $m$ large. 
\item[{\rm{(2)}}] If $j_b(z_0) =0$, then $p_{mp+b-1}(z)$ has a zero {\rm{(}}$k$ zeros 
if $z$ has a $k$-th order zero at $z_0${\rm{)}} near $z_0$ for $m$ large. 
\item[{\rm{(3)}}] There are at most $2p+2b-3$ points in $\bbC\backslash\text{bands}$ 
where $j_b (z_0)$ is zero. 
\end{SL} 
\end{theorem} 

\begin{proof} \eqref{2.39} is immediate from \eqref{2.34} and $\abs{\Gamma_-/\Gamma_+}<1$. 
(1) and (2) then follow by Hurwitz's theorem if we show that $j_b(z)$ is not identically 
zero. 

By \eqref{2.1} and \eqref{2.7} near $z=\infty$, 
\[ 
T_p(z) = \biggl( \, \prod_{j=1}^p a_j \biggr)^{-1} z^p 
\begin{pmatrix} 1 & 0 \\ 0 & 0 \end{pmatrix} +O(z^{p-1}) 
\]
which implies $\Gamma_+ = (\prod_{j=1}^p a_j)^{-1} z^p + O(z^{p-1})$ and $\Gamma_-(z) 
=O(z^{-p})$. It follows that $a(z)\to 0$ as $z\to\infty$ and $b(z)\to 1$. Thus, since 
$p_{b-1}$ has degree $b-1$, \eqref{2.34} shows that as $z\to\infty$ on the main sheet, 
$f(z_0)$ has a pole of order $b-1$. 

On the other sheet, $P_+$ changes to $P_-$, so $a(z)\to 0$ and $b(z)\to 0$ on the other 
sheet. It follows that $j(z)$ has a pole at $\infty$ of degree at most $b-2$. $j$ 
also has poles of degree at most $1$ at each branch point. Thus, $j_b(z)$ as a 
function on $\calS$ has total degree at most $2p+ (b-1) + (b-2) =2p + 2b-3$ which 
bounds the number of zeros. 
\end{proof} 

Finally, we turn to zeros on the bands. A major role will be played by the 
function on the right side of \eqref{2.39} ($j$ is for ``Jost" since 
this acts in many ways like a Jost function): 
\begin{equation} \lb{2.40} 
j_b(z) = a(z) q_{b-2}(z) + b(z) p_{b-1}(z) 
\end{equation} 

\begin{lemma} \lb{L2.4} $j_b$ is nonvanishing on the interior of the bands. 
\end{lemma} 

{\it Remark.} By $j_b(x)$ for $x$ real, we mean \eqref{2.40} with $a$ defined via 
$\lim_{\veps\downarrow 0} a(x+i\veps)$ since $P_\pm$ are only defined off 
$\bbC\backslash\text{bands}$. 

\begin{proof} As already mentioned, the boundary values obey 
\begin{equation} \lb{2.41} 
\lim_{\veps\downarrow 0}\, P_+ (x+i\veps) = \lim_{\veps\downarrow 0}\, P_- (x-i\veps) 
\end{equation} 
(by the two-sheeted nature of $P_+$ and $P_-$). Thus, by \eqref{2.26} and \eqref{2.27}, 
\begin{align} 
a(x+i0) &= -a (x-i0) \lb{2.42} \\
b(x+i0) &= 1-b(x-i0) \lb{2.43} 
\end{align}
Moreover, since $T_p$ and $\Gamma_\pm$ are real on $\bbR\backslash\text{bands}$, 
$a(z)$ and $b(z)$ are real on $\bbR\backslash\text{bands}$ (by \eqref{2.22}). Thus 
\begin{align}
a(x+i0) &= \ol{a(x-i0)}  \lb{2.44} \\
b(x+i0) &= \ol{b(x-i0)} \lb{2.45} 
\end{align} 
The last four equations imply for $x$ in the bands
\begin{align} 
\Real(a(x+i0)) &= 0 \lb{2.46} \\ 
\Real (b(x+i0)) &= \tfrac12  \lb{2.47} 
\end{align} 

$p$ and $q$ are real on $\bbR$, so 
\begin{equation} \lb{2.48} 
\Real (j_b(x)) = \tfrac12\,  p_{b-1}(x)
\end{equation} 
Thus, if $j_b(x_0) =0$ on the bands, $p_{b-1}(x_0) =0$. 

As we have seen, $a(z)=0$ only at the Dirichlet points and so not in the bands. 
If $p_{b-1}(x_0)=0=j_b(x_0)$, then since $a(x_0)\neq 0$, we also have $q_{b-2}(x_0) 
=0$. By \eqref{2.10a}, if $p_{b-1}(x_0) =q_{b-2}(x_0)$, then $\det (T_b (x_0))=0$, 
which is false. We conclude via proof by contradiction that $j_b(x)$ has no zeros. 
\end{proof} 

\begin{theorem}\lb{T2.5} For each $b$ and each band $j$, there is an integer 
$D_{b,j}$ so the number of zeros $N_{b,j}(m)$ of $p_{mp+b-1}$ is either 
$m-D_{b,j}$ or $m-D_{b,j}+1$. In particular, 
\begin{equation} \lb{2.49} 
\sup_{n,j}\, \biggl| \f{n}{p} - N^{(n,j)}\biggr| < \infty 
\end{equation} 
Moreover, \eqref{1.8} holds. 
\end{theorem} 

\begin{proof} By \eqref{2.34}, \eqref{2.41}, \eqref{2.42}, and \eqref{2.43}, we have 
\begin{equation} \lb{2.50} 
p_{mp+b-1}(x) = j_b(x) \Gamma_+(x)^m + \ol{j_b(x)}\,\, \ol{\Gamma_+(x)}^m 
\end{equation} 
on the bands. By the lemma, $j_b(x)$ is nonvanishing inside band $j$, so 
\begin{equation} \lb{2.51} 
j_b(x) = \abs{j_b(x)} e^{i\gamma_b(x)} 
\end{equation} 
where $\gamma_b$ is continuous --- indeed, real analytic --- and by a simple argument, 
$\gamma_b$ and $\gamma'_b$ have limits as $x\downarrow a_j$ or $x\uparrow b_j$. 

By \eqref{2.37}, \eqref{2.50} becomes 
\begin{equation} \lb{2.52} 
p_{mp+b-1}(x) = 2\abs{j_b(x)} \cos (\pi mp\, k(x) + \gamma_b(x)) 
\end{equation} 
Define $D_{b,j}$ to be the negative of the integral part of $[\gamma_b(b_j)- 
\gamma_b(a_j)]/\pi$. Since $\sup_{\text{bands}} \abs{\gamma'_b(x)}<\infty$, there 
is, for large $m$, at most one solution of $\pi mp\, k(x) + \gamma_b(x) =\pi\ell$ 
for each $\ell$. Given this, it is immediate that the number of zeros is $m-D_{b,j}$ 
or $m-D_{b,j}+1$. 

Finally, \eqref{1.8} is immediate from \eqref{2.52}. Given that $\gamma$ is $C^1$, 
we even get that 
\begin{equation} \lb{2.53} 
k(x_{\ell+1}^{(n,j)}) - k(x_\ell^{(n,j)}) = \f{1}{n} + O\biggl( \f{1}{n^2}\biggr) 
\end{equation} 
\end{proof} 

As for point (3) from the introduction, the proof of Theorem~\ref{T2.4} shows that if 
$z_0$ is not in the bands and is a limit of zeros of $p_{mp+b-1}(z)$, then $p_{mp+b-1}
(z_0)$ goes to zero exponentially (as $\Gamma_-^m$). If this is true for each $b$, 
then $\sum_{n=0}^\infty \abs{p_n(z)}^2 <\infty$, which means $z_0$ is in the pure 
point spectrum of $d\mu$. Since the bands are also in the spectrum, we have 

\begin{proposition} \lb{P2.6} $z_0\in\bbC$ is a limit of zeros of $p_n(z)$ 
{\rm{(}}all $n${\rm{)}} if and only if $z_0\in\supp(d\mu)$. 
\end{proposition} 

{\it Remark.} This also follows from a result of Denisov-Simon \cite{DS2}, but their 
argument, which applies more generally, is more subtle.

\section{OPUC With Periodic Verblunsky Coefficients} \lb{s3} 

In this section, we analyze the zeros of OPUC with Verblunsky coefficients obeying 
\eqref{1.6}. We begin with a summary of the transfer matrices, discriminants, and
Abelian functions in this situation. These ideas, while an obvious analog of the 
OPRL situation, seem not to have been studied before their appearance in \cite{OPUC2},  
which is the reference for more details. Many of the consequences of these ideas 
were found earlier in work of Peherstorfer and Steinbauer 
\cite{PS6,PS1,PS2,PS3,PS4,PS2000,PS00}. 

Throughout, we will suppose that $p$ is even. If $(\alpha_0, \dots, \alpha_{p-1}, 
\alpha_p, \dots)$ is a sequence with odd period, $(\beta_0, \beta_1, \dots) = 
(\alpha_0, 0, \alpha_1, 0, \alpha_2, \dots)$ has even period and 
\begin{equation} \lb{3.1} 
\Phi_{2n} (z,\{\beta_j\}) = \Phi_n (z^2, \{\alpha_j\}) 
\end{equation} 
so results for the even $p$ case immediately imply results for the odd $p$. 

Define the $2\times 2$ matrix 
\begin{equation} \lb{3.2} 
A_k(z) = \f{1}{\rho_k} \begin{pmatrix} 
z & -\bar\alpha_k \\ -z\alpha_k & 1 \end{pmatrix} 
\end{equation} 
where $\rho_k$ is given by \eqref{1.4}. Then 
\[
\det (A_k (\alpha)) = z 
\]
\eqref{1.2} and its $^*$ are equivalent to 
\begin{equation} \lb{3.3} 
\binom{\varphi_{n+1}}{\varphi_{n+1}^*} = A_n (z) \binom{\varphi_n}{\varphi_n^*} 
\end{equation} 
The second kind polynomials, $\psi_n(z)$, are the OPUC with Verblunsky coefficients 
$\{-\alpha_j\}_{j=0}^\infty$. Then it is easy to see that 
\begin{equation} \lb{3.4} 
\binom{\psi_{n+1}}{-\psi_{n+1}^*} = A_n(z) \binom{\psi_n}{-\psi_n^*}  
\end{equation} 
with $A$ given by \eqref{3.2}. 

We thus define 
\begin{equation} \lb{3.5} 
T_n(z) = A_{n-1}(z) \dots A_0(z) 
\end{equation} 
By \eqref{1.6}, we have 
\begin{equation} \lb{3.6} 
T_{mp+b} = T_b (T_p)^m 
\end{equation}  
\eqref{3.3} and \eqref{3.4} imply that 
\begin{align} 
\binom{\varphi_n}{\varphi_n^*} &= T_n \binom{1}{1} \lb{3.7} \\
\binom{\psi_n}{-\psi_n^*} &= T_n \binom{1}{-1} \lb{3.8} 
\end{align} 
so that 
\begin{equation} \lb{3.9} 
T_n(z) = \tfrac12 \begin{pmatrix} 
\varphi_n(z) + \psi_n(z) & \varphi_n(z) -\psi_n(z) \\
\varphi_n^*(z) - \psi_n^*(z) & \varphi_n^*(z) + \psi_n^*(z) 
\end{pmatrix} 
\end{equation} 

The discriminant is defined by 
\begin{equation} \lb{3.10} 
\Delta(z) = z^{-p/2} \, \tr (T_p(z)) 
\end{equation} 
The $z^{-p/2}$ factor (recall $p$ is even) is there because $\det (z^{-p/2} T_p(z)) 
=1$, so $z^{-p/2} T_p(z)$ has eigenvalues $\Gamma_\pm (z)$ given by \eqref{2.11}. 
$\Delta(z)$ is real on $\partial\bbD$ so 
\begin{equation} \lb{3.10a} 
\Delta (z) = \ol{\Delta(1/\bar z)}
\end{equation} 

$\Delta(z)\in (-2,2)$ only if $z=e^{i\theta}$ and there are $p$ roots, each of 
$\tr (T_p(z)) \mp 2z^{p/2}=0$, that is, $p$ solutions of $\Delta(z) =\pm 2$. These 
alternate on the circle at points $+2, -2, -2, +2, +2, -2, -2, \dots$, so we pick 
\begin{equation} \lb{3.11} 
0\leq x_1 < y_1 \leq x_2 < y_2 \leq \cdots < y_p \leq 2\pi
\end{equation} 
where $e^{ix_j}, e^{iy_j}$ are solutions of $\Delta(z) =\pm 2$. 

The bands 
\begin{equation} \lb{3.12} 
B_j =\{e^{i\theta}\mid x_j \leq \theta \leq y_j\}
\end{equation} 
are precisely the points where $\Delta(z)\in [-2,2]$. In between are the gaps 
\begin{equation} \lb{3.13} 
G_j \{e^{i\theta}\mid y_j < \theta < x_{j+1}\}
\end{equation} 
where $x_{p+1} = x_1 + 2\pi$. Some gaps can be closed, that is, $G_j$ is empty 
(i.e., $y_j =x_{j+1}$). 

We also see that on $\bbC\backslash\text{bands}$, $\abs{\Gamma_+} > \abs{\Gamma_-}$, 
so the Lyapunov exponent is given by 
\begin{equation} \lb{3.14} 
\lim_{n\to\infty}\, \tfrac{1}{n}\, \log \|T_n(z)\| = \tfrac12 \, \log \abs{z} + 
\tfrac{1}{p}\, \log \abs{\Gamma_+(z)} \equiv \gamma(z)
\end{equation} 

If we remove the bands from $\bbC$, \eqref{2.12x} holds. Moreover, $\Gamma_+(z)$ has 
an analytic continuation to the Riemann surface, $\calS$, of $[\prod_{\text{open gaps}}
(z-e^{iy_{j+1}})(z-e^{ix_j})]^{1/2}$. The genus of $\calS$, $\ell\leq p-1$, where 
$\ell+ 1$ is the number of open gaps. (In some sense, the OPRL case, where the genus 
$\ell$ is the number of gaps, has $\ell +1$ gaps also, but one gap is $\bbR\backslash 
[\alpha_1, \beta_p]$ which includes infinity.) $\Gamma_-$ is the analytic continuation 
of $\Gamma_+$ to the second sheet. 

The Dirichlet data are partly these points in $\partial\bbD$, $z_j$, 
\begin{equation} \lb{3.15} 
T_p(z) \binom{1}{1} = c_z \binom{1}{1}
\end{equation} 
It can be shown there is one such $z_j$ in each gap (including closed gaps) for 
the $p$ roots of $\varphi_p(z) - \varphi_p^*(z)$. We let $c_j = c_{z_j}$. If $z_j$ 
is at a gap edge, $\abs{c_j}=1$; otherwise $\abs{c_j}\neq 1$. If $\abs{c_j} >1$, 
we add sign $-1$ to $z_j$ and place the Dirichlet point on the lower sheet of $\calS$ 
at point $z_j$. If $\abs{c_j}<1$, we add sign $+ 1$ and put the Dirichlet point 
on the initial sheet. $+1$ points correspond to pure points in $d\mu$. 

As in the OPRL case, the set of possible Dirichlet data points is a torus, but 
now of dimension $\ell +1$. This torus parametrizes those $\mu$ with periodic 
$\alpha$'s and discriminant $\Delta$. 

The $F$-function, \eqref{1.10}, has a meromorphic contribution to $\calS$ with poles 
precisely at the Dirichlet data points. 

The potential theoretic equilibrium measures $dk$ for the bands have several critical 
properties: 
\begin{SL} 
\item[(1)] If $k(e^{i\theta_0}) = k(\{e^{i\theta}\mid x_1 < \theta < \theta_0\})$, then 
\begin{equation} \lb{3.16} 
k(e^{iy_j}) = k(e^{ix_{j+1}}) = \f{j}{p}
\end{equation} 
\item[(2)] The Thouless formula holds: 
\begin{equation} \lb{3.17} 
\gamma(z) = \int \log \abs{z-e^{i\theta}}\, dk(e^{i\theta}) + \log C_B
\end{equation} 
where $\gamma$ is given by \eqref{3.14} and $C_B$ is the capacity of the bands. 
\item[(3)] We have 
\begin{equation} \lb{3.18} 
C_B =\prod_{j=0}^{p-1} (1-\abs{\alpha_j}^2)^{1/2}
\end{equation} 
\item[(4)] 
\begin{equation} \lb{3.19} 
\Gamma_+ (z) = C_B z^{-p/2} \exp \biggl( p\int \log (z-e^{i\theta})\, dk (e^{i\theta})\biggr)
\end{equation} 
\end{SL}

This completes the review of periodic OPUC. The analog of Theorem~\ref{T2.1} does not 
involve $\Phi_n$ but $\Phi_n - \Phi_n^*$:  

\begin{theorem} \lb{T3.1} The zeros of $\Phi_{mp}(z) - \Phi_{mp}^*(z)$ are at the 
following points: 
\begin{SL} 
\item[{\rm{(i)}}] the $p$ Dirichlet data $z_j$'s in each gap of the period $p$ problem. 
\item[{\rm{(ii)}}] the $(m-1)p$ points where 
\begin{equation} \lb{3.21} 
k(e^{i\theta}) = \f{k-1}{p} + \f{q}{mp} 
\end{equation} 
$k=1, \dots, p$; $q=1, \dots, m-1$. 
\end{SL} 
\end{theorem} 

\begin{proof} As noted (and proven several ways in \cite[Chapter~11]{OPUC2}), for a period 
$mp$ problem, $\Phi_{mp}-\Phi_{mp}^*$ has its zeros, one in each gap. The gaps of the 
$mp$ problem are the gaps of the original problem plus a closed gap at each point where 
\eqref{3.21} holds. There is a zero in each closed gap and at each point where \eqref{3.15} 
holds since then $T_{mp}(z) \binom{1}{1} = c_j^m \binom{1}{1}$. 
\end{proof} 

We now turn to the analysis of zeros of $\varphi_{mp+b}(z)$, $b=0,1,\dots, p-1$; 
$m=0,1,2,\dots$. The analog of \eqref{2.37a} is, by \eqref{3.7}, 
\begin{equation} \lb{3.22} 
\varphi_{mp+b} =\bigg\langle \binom{1}{0}, T_b (T_p)^m \binom{1}{1}\bigg\rangle 
\end{equation} 
As in Section~\ref{s2}, we write, for $z\in\bbC\backslash\text{bands}$: 
\begin{equation} \lb{3.23} 
z^{-p/2} T_p(z) = \Gamma_+(z) P_+(z) + \Gamma_-(z) P_-(z) 
\end{equation} 
where $P_\pm$ are $2\times 2$ matrices which are complementary projections, that is, 
\eqref{2.24}/\eqref{2.25} hold. \eqref{2.26}/\eqref{2.27} are replaced by 
\begin{align} 
P_+ &= \f{z^{-p/2} T_p(z) - \Gamma_-(z)\boldsymbol{1}}{\Gamma_+ - \Gamma_-} \lb{3.24} \\
P_- &= \f{z^{-p/2} T_p(z) - \Gamma_+(z)\boldsymbol{1}}{\Gamma_- - \Gamma_+} \lb{3.25} 
\end{align} 
So, in particular, $P_\pm$ have meromorphic continuations to $\calS$, and $P_+$ 
continued to the other sheet is $P_-$. 

Define 
\begin{align} 
a(z) &= \f12\, \bigg\langle\binom{1}{1}, P_+ \binom{1}{1}\bigg\rangle \lb{3.26} \\ 
b(z) &= \f12\, \bigg\langle\binom{1}{-1}, P_+ \binom{1}{1}\bigg\rangle \lb{3.27} 
\end{align}
so that, by \eqref{2.25}, 
\begin{align} 
\f12\, \bigg\langle\binom{1}{1}, P_-\binom{1}{1}\bigg\rangle &= 1-a(z) \lb{3.28} \\ 
\f12\, \bigg\langle\binom{1}{-1}, P_-\binom{1}{1}\bigg\rangle &= -b(z) \lb{3.29}
\end{align} 
Thus, since $\f{1}{\sqrt2}\binom{1}{1}, \f{1}{\sqrt2}\binom{1}{-1}$ are an orthonormal 
basis, 
\begin{equation} \lb{3.30} 
\begin{split}
z^{-mp/2} T_{mp}\binom{1}{1} &= \Gamma_+^m \biggl[ a(z) \binom{1}{1} + b(z) \binom{1}{-1}\biggr] \\
& \qquad \qquad + \Gamma_-^m \biggl[ (1-a(z))\binom{1}{1} - b(z)\binom{1}{-1}\biggr] 
\end{split}
\end{equation} 
Therefore, by \eqref{3.7}, \eqref{3.8}, and \eqref{3.22}, 
\begin{equation} \lb{3.31} 
\begin{split} 
\varphi_{mp+b}(z) &= \varphi_b(z) [a(z)z^{mp/2}\Gamma_+^m + (1-a) z^{mp/2}\Gamma_-^m] \\
& \qquad \qquad + \psi_b(z) [b(z) z^{mp/2} \Gamma_+^m -b(z) z^{mp/2}\Gamma_-^m] 
\end{split}
\end{equation} 
We thus define 
\begin{equation} \lb{3.32} 
j_b(z) =a(z) \varphi_b(z) + b(z) \psi_b(z) 
\end{equation} 
and we have, since $\abs{\Gamma_+} > \abs{\Gamma_-}$ on $\bbC\backslash\text{bands}$: 

\begin{theorem} \lb{T3.2} For $z\in\bbC\backslash\text{bands}$, 
\begin{equation} \lb{3.33} 
\lim_{m\to\infty}\, z^{-mp/2} \Gamma_+^{-m} \varphi_{mp+b}(z) = j_b(z) 
\end{equation} 
In addition, $j_b$ is nonvanishing near $z=\infty$. 

In particular, if $z_0\notin\text{bands}$ and $j_b (z_0)\neq 0$, then for some $\veps >0$  
and $M$\!, we have $\varphi_{mp+b}(z_0)\neq 0$ if $\abs{z-z_0}<\veps$ and $m\geq M$\!. 
If $z_0\notin\text{bands}$ and $j_b(z_0)$ has a zero of order $k$, then for some 
$\veps >0$ and all $m$ large, $\varphi_{mp+b}(z)$ has precisely $k$ zeros {\rm{(}}counting  
multiplicity{\rm{)}}. The number of $z_0$ in $\bbC\backslash\text{bands}$ with $j_b(z_0)=0$ 
is at most $2p + 2b-1$. 
\end{theorem} 

\begin{proof} As noted, \eqref{3.31} and $\abs{\Gamma_+} > \abs{\Gamma_-}$ imply 
\eqref{3.33}. To analyze $j_b(z)$ near $z=\infty$, we proceed as follows: We have, by 
\eqref{3.2} and \eqref{3.5}, that as $\abs{z}\to\infty$, 
\begin{align} 
T_p(z) &= z^p \biggl( \, \prod_{j=0}^{p-1} \rho_j^{-1}\biggr) 
\biggl[ \begin{pmatrix} 1 & 0 \\ -\alpha_{p-1} & 0 \end{pmatrix} \dots 
\begin{pmatrix} 1 & 0 \\ -\alpha_0 & 0 \end{pmatrix} \biggr] + O(z^{p-1}) \lb{3.34} \\
&= z^p \biggl(\, \prod_{j=1}^{p-1} \rho_j^{-1} \biggr) 
\begin{pmatrix} 1 & 0 \\ -\alpha_{p-1} & 0 \end{pmatrix} + O(z^{p-1}) \lb{3.35} 
\end{align} 
from which it follows that 
\begin{align} 
P_+ &= \begin{pmatrix} 1 & 0 \\ -\alpha_{p-1} & 0 \end{pmatrix} + O(z^{-1}) \lb{3.36} \\ 
P_- &= \begin{pmatrix} 0 & 0 \\ \alpha_{p-1} & 1 \end{pmatrix} + O(z^{-1}) \lb{3.37} 
\end{align} 
and 
\begin{align} 
a(z) &= \tfrac12\, (1-\alpha_{p-1}) + O(z^{-1}) \lb{3.38} \\
b(z) &= \tfrac12\, (1+\alpha_{p-1}) + O(z^{-1}) \lb{3.39} 
\end{align} 

We have 
\begin{align} 
\varphi_b(z) &= \biggl(\, \prod_{j=0}^{p-1} \rho_j^{-1}\biggr) z^b + O(z^{b-1}) \lb{3.40} \\ 
\psi_b(z) &= \biggl(\, \prod_{j=0}^{p-1} \rho_j^{-1}\biggr) z^b + O(z^{b-1}) \lb{3.41}  
\end{align} 
from which we see that 
\begin{equation} \lb{3.42} 
j_b(z) = \biggl(\, \prod_{j=0}^{p-1} \rho_j^{-1}\biggr) z^b + O(z^{b-1})  
\end{equation} 
since \eqref{3.38}/\eqref{3.39} imply $a(z) + b(z) =1+O(z^{-b})$. In particular, 
$j_b(z)$ is not zero near $\infty$, so $j_b$ is not identically zero, and the assertion 
about locations of zeros of $\varphi_{mp+b}(z)$ follows from Hurwitz's theorem. 

Since $1-a(z) -b(z)=O(z^{-1})$, \eqref{3.28}/\eqref{3.29} imply that, on the second sheet,  
the analytic continuation of $j_b(z)$ near $\infty$ is $O(z^{b-1})$. It follows that 
$j_b$ has a pole of order $b$ at $\infty$ on the main sheet (regular if $b=0$) and a 
pole of order at most $b-1$ (a zero if $b=0$ and is regular if $b=1$) at $\infty$ 
on the second sheet. $j_b$ also can have at most $2p$ simple poles at the $2p$ branch 
points. 

It follows that the degree of $j_b$ as a meromorphic function on $\calS$ is at most 
$2p+2b -1$ (if $b=0$, $2p$). Thus the number of zeros is at most $2p+2b-1$ if 
$b\neq 0$. If $b=0$, there are at most $2p + 2b$ zeros. But since then one is 
at $\infty$ on the second sheet, the number of zeros on finite points is at most 
$2p+2b -1$. 
\end{proof} 

Next, we note that 

\begin{theorem} \lb{T3.3} Let $\{\alpha_n\}$ be periodic and not at all $0$. $z_0$ is a 
limit of zeros of $\varphi_n(z)$ {\rm{(}}i.e., there exist $z_n$ with $\varphi_n (z_n) 
=0$ and $z_n \to z_0$ if and only if $z_0$ lies in the support of $d\mu${\rm{)}}. 
\end{theorem} 

{\it Remark.} $\alpha_n =0$ has $0$ as a limit point of zeros at $\varphi_n(z) 
=z^n$, so one needs some additional condition on the $\alpha$'s to assure this result. 

\begin{proof} By Theorems~8.1.11 and 8.1.12 of \cite{OPUC1}, if $z_0\in\supp(d\mu)$, 
then it is a limit point of zeros. For the other direction, suppose $z_0\notin\text{bands}$ 
and is a limit point of zeros. By Theorem~\ref{T3.2}, $j_b(z_0)=0$ for each $b=0,1, 
\dots, p-1$, so by \eqref{3.31}, $\varphi_{mp+b}(z)\sim C(\Gamma_- z_0^{p/2})^m$ 
which, since $\abs{z_0}\leq 1$ and $\abs{\Gamma_-}<1$, implies that $\varphi_n (z_0)$ 
goes to zero exponentially. 

Since $\alpha_n$ is not identically zero, some $\alpha_j$, $j\in \{0,1,\dots, p-1\}$ 
is nonzero. Thus, by Szeg\H{o} recursion for $\varphi_j$, 
\[
\varphi_{mp+j}^*(z_0) = \alpha_j^{-1} [z_0 \varphi_{mp}(z_0) - \rho_j 
\varphi_{mp+1}(z_0)] 
\]
goes to zero exponentially in $m$. 

Since $\alpha_n$ is periodic, $\sup_n \abs{\alpha_n}<1$, and so, $\sup_n \rho_n^{-1} 
<\infty$. Since 
\[
\varphi_{mp+j+1}^*(z_0) =\rho_{j+1}^{-1} (\varphi_{mp+j}^*(z_0) - 
\alpha_j \varphi_{mp+j}(z_0)) 
\]
we see $\varphi_{mp+j+1}^*(z_0)$ decays exponentially and so, by induction, 
$\varphi_n^*(z_0)$ decays exponentially. By the Christoffel-Darboux formula 
(see \cite[eqn.~(2.2.70)]{OPUC1}), $\abs{\varphi_n^*(z_0)}^2 \geq 1-\abs{z_0}^2$, 
so the decay implies $\abs{z_0}=1$. But if $z_0\in\partial\bbD$ and $\sum_n 
\abs{\varphi_n(z_0)}^2 <\infty$, then $\mu (\{z_0\})>0$ 
(see \cite[Theorem~2.7.3]{OPUC1}). 

Thus if $z_0$ is a limit of zeros, either $z_0\in\text{bands}$ or $\mu (\{z_0\})>0$, that is, 
$z_0\in\supp (d\mu)$. \
\end{proof} 

Finally, in our analysis of periodic OPUC, we turn to zeros near to the bands. We define 
$\ti\jmath_b$ on $\bbC\backslash\text{bands}$ so that \eqref{3.31} becomes 
\begin{equation} \lb{3.43} 
\varphi_{mp+b}(z) = j_b(z) z^{mp/2} \Gamma_+^m + \ti\jmath_b(z) z^{mp/2} \Gamma_-^m 
\end{equation} 
While $\varphi_{mp+b}(z)$ is continuous across the bands, $j_b$, $\ti\jmath_b$, and 
$\Gamma_\pm$ are not. In fact, $\Gamma_+$ (resp.~$j_b$) continued across a band 
becomes $\Gamma_-$ (resp.~$\ti\jmath_b$). We define all four objects at $e^{i\theta}
\in\partial\bbD$ as limits as $r\uparrow 1$ of the values at $re^{i\theta}$.  

\begin{proposition} \lb{P3.4} 
\begin{SL} 
\item[{\rm{(i)}}] In the bands, 
\begin{equation} \lb{3.44} 
e^{i\theta p/2}\Gamma_+ (e^{i\theta}) = \exp (-i\pi p \, k(\theta)) 
\end{equation} 
\item[{\rm{(ii)}}] At no point in the bands do both $j_b(e^{i\theta})$ and $\ti\jmath_b 
(e^{i\theta})$ vanish. 
\item[{\rm{(iii)}}] $\ti\jmath_b$ is everywhere nonvanishing on the interiors 
of the bands. 
\end{SL} 
\end{proposition} 

\begin{proof} (i) This follows from \eqref{3.19}. There is an issue of checking that 
it is $\exp (-i\pi p\, k(\theta))$, not $\exp (i\pi p\, k(\theta))$. To confirm this, 
note that $\f{\partial}{\partial\theta} \Ima \log (\exp(-i\pi p \, k(\theta)))\leq 0$ 
and mainly $<0$. Since $\partial\abs{\Gamma_+}/\partial r \leq 0$ at $r=1$, this is 
consistent with \eqref{3.44} and the Cauchy-Riemann equations. 

\smallskip 
(ii) follows from \eqref{3.43} and the fact that $\varphi_n(z)$ is nonvanishing on 
$\partial\bbD$.  

\smallskip 
(iii) Continue \eqref{3.43} through the cut. Since $\varphi_m$ is entire, the 
continuation onto the ``second sheet" is also $\varphi_m$. $\Gamma_\pm$ get 
interchanged by crossing the cut. Let us use $j_{b,2},\ti\jmath_{b,2}$ for the 
continuation to the second sheet (of course, $j_{b,2}$ is $\ti\jmath_b$ on the 
second sheet, but that will not concern us). 

By this \eqref{3.43} continued, $\varphi_{mp+b}(z)=0$ if and only if 
\begin{equation} \lb{3.45} 
\biggl( \f{\Gamma_-(z)}{\Gamma_+(z)}\biggr)^m = - \f{\ti\jmath_{b,2}(z)}{j_{b,2}(z)}  
\end{equation} 

If $\ti\jmath_b(z_0) =0$ for $z_0\in\partial\bbD$, then 
$\abs{\ti\jmath_{b,2}(rz_0)/j_{b,2}(rz_0)}$ goes from $0$ to a nonzero 
value as $r$ increases. On the other hand, since $\abs{\Gamma_-/\Gamma_+} 
<1$ on $\bbC\backslash\text{bands}$, for $m$ large, 
$\abs{\Gamma_-(rz_0)/\Gamma_+ (rz_0)}^m$ goes from $1$ to a very 
small value as $r$ increases. It follows that for $m$ large, 
\[ 
\biggl| \f{\Gamma_-(z)}{\Gamma_+ (z)}\biggr|^m = 
\biggl| \f{\ti\jmath_{b,2}(z)}{j_{b,2}(z)}\biggr| 
\]
has a solution $r_m z_0$ with $r_m >1$ and $r_m \to 1$. As in \cite{SaffProc}, 
we can change the phase slightly to ensure \eqref{3.45} holds for some point, 
$z_m$, near $r_m z_0$ with $\abs{z_m} >1$. Since $\varphi$ has no zero in 
$\bbC\backslash\bbD$, this is a contradiction. 
\end{proof}  

{\it Remark.} This proof shows that in the bands $\abs{\ti\jmath_b (e^{i\theta})} > 
\abs{j_b (e^{i\theta})}$. 

\smallskip
\eqref{3.43} says we want to solve 
\begin{equation} \lb{3.48} 
\biggl( \f{\Gamma_-(z)}{\Gamma_+(z)}\biggr)^m = g_b(z) 
\end{equation}  
to find zeros of $\varphi_{mp+b}(z)$. We have, by the remark, that 
$\abs{g(\theta)} <1$. 

\smallskip
\noindent{\bf Definition.} We call $z_0\in\text{bands}$ a singular point of order 
$k$ if $j_b(z_0)=0$ and the zero is of order $k$.  

\smallskip 
We do not know if there are singular points in any example! If so, they should 
be nongeneric. We define the functions 
\begin{equation} \lb{3.46} 
\ti g_b(\theta) = -\f{j_b (e^{i\theta})}{\ti\jmath_b (e^{i\theta})}
\end{equation} 
and 
\begin{equation} \lb{3.47} 
g_b(z) = -\f{j_b(z)}{\ti\jmath_b(z)} 
\end{equation}  
For $e^{i\theta}$ in the interior of a band minus the singular points, let 
$A(\theta)$ be given by 
\begin{equation} \lb{3.49} 
\f{\ti g_b(\theta)}{\, \ol{\ti g_b (\theta)}\,} = \exp (2i A(\theta)) 
\end{equation} 
with $A$ continuous away from the singular points. 

The analysis of a similar equation to \eqref{3.48} in \cite{SaffProc} shows that:  
\begin{SL} 
\item[(a)] The solutions of \eqref{3.48} near $\abs{z}=1$ lie in sectors 
where 
\begin{equation} \lb{3.50} 
2\pi p \, k(\theta) = A(\theta) + \f{2\pi j}{m} + O\biggl( \f{1}{m\log m}\biggr) 
\end{equation}  
with exactly one solution in each such sector. 

\item[(b)] The magnitudes of the solutions obey 
\begin{equation} \lb{3.51} 
\abs{z} = 1 - O\biggl( \f{\log m}{m}\biggr)  
\end{equation} 

\item[(c)] Successive zeros $z_{k+1},z_k$ obey  
\begin{equation} \lb{3.52} 
k (\arg (z_{k+1})) - k(\arg z_k) = \f{1}{mp} + O\biggl( \f{1}{m\log m}\biggr) 
\end{equation} 
and  
\begin{equation} \lb{3.53} 
\biggl| \f{z_{k+1}}{z_k}\biggr| = 1 + O\biggl( \f{1}{m\log m}\biggr) 
\end{equation} 

\item[(d)] All estimates in (a)--(c) are uniform on a  band. 

\item[(e)] Away from singular points, all $O(1/m\log m)$ errors can be replaced by 
$O(1/m^2)$ and $O(\log m/m)$ in \eqref{3.51} by $O(1/m)$. If there are no singular 
points, these are uniform over a band. 
\end{SL} 

\smallskip
It is easy to see that the total variation of $A$ in each interval between 
singular points (or band endpoints) is finite, so \eqref{3.50} and the fact 
that $k$ varies by $1/p$ over a band say that the number of solutions in a band 
differs from $m$ by a finite amount. This implies 

\begin{theorem}\lb{T3.5} Let $N^{(n,j)}$ be the number of zeros, $z_0$, of 
$\varphi_n(z)$ that obey 
\begin{SL} 
\item[{\rm{(a)}}] $\arg z_0 \in \text{band $j$}$ 
\item[{\rm{(b)}}] 
\begin{equation} \lb{3.54} 
(1-\abs{z_0}) \leq n^{-1/2} 
\end{equation} 
\end{SL} 
Then 
\begin{SL} 
\item[{\rm{(a)}}] $\sup_{n,j} \abs{N^{(n,j)} - \f{n}{p}} < \infty$ 
\item[{\rm{(b)}}] For $n$ large, all such zeros have 
\begin{equation} \lb{3.55} 
(1-\abs{z_0}) \leq C \, \f{\log n}{n}  
\end{equation} 
and if there are no singular points, we can replace $\log n/n$ in \eqref{3.55} by 
$1/n$. 
\end{SL} 
\end{theorem}

\bigskip

\end{document}